\documentclass[10pt]{amsart}
\usepackage{latexsym}
\usepackage[dvips]{color}
\usepackage{float}
\usepackage{afterpage} 
\usepackage [pdftex]{graphicx}
\usepackage{epsfig}
\usepackage{amsmath}
\usepackage{amsthm}
\usepackage{amsfonts}
\usepackage{amssymb}
\usepackage{setspace}
\usepackage{subfigure}
\usepackage{tikz}
\usetikzlibrary{arrows}
\usetikzlibrary{arrows.meta}
\usetikzlibrary{matrix}
\usepackage{subfigure}

\newtheorem{theorem}{Theorem}
\newtheorem{definition}{Definition}
\newtheorem{proposition}{Proposition}
\newtheorem{example}{Example}

\newtheorem{corollary}{Corollary}

\begin{document}

\title[Hidden categories: Lewin's GIS and K-nets]{Hidden categories: a new perspective on Lewin's Generalized Interval Systems and Klumpenhouwer networks.}
\author{Alexandre Popoff}
\address{Independent Researcher}
\email{al.popoff@free.fr}

\author{Moreno Andreatta}
\address{CNRS/Institute for Advanced Mathematical Research, ITI CREAA, University of Strasbourg, France and IRCAM, Paris, France}
\email{andreatta@math.unistra.fr, andreatta@ircam.fr}


\begin{abstract}
In this work we provide a categorical formalization of several constructions found in transformational music theory. We first revisit David Lewin's original theoretical construction of Generalized Interval Systems (GIS) to show that it implicitly defines categories. When all the conditions in Lewin's definition are fullfilled, such categories coincide with the category of elements $\int_\mathbf{G} S$ for the group action $S \colon \mathbf{G} \to \mathbf{Sets}$ implied by the GIS structure. By focusing on the role played by categories of elements in such a context, we reformulate previous definitions of transformational networks in a $\mathbf{Cat}$-based diagrammatical perspective, and present a new definition of transformational networks (called \textit{CT-Nets}) in general musical categories. We show incidently how such an approach provides a bridge between algebraic, geometrical and graph-theoretical approaches in transformational music analysis. We end with a discussion on the new perspectives opened by such a formalization of transformational theory, in particular with respect to $\mathbf{Rel}$-based transformational networks which occur in well-known music-theoretical constructions such as Douthett's and Steinbach's Cube Dance.
\end{abstract}

\maketitle
\section{Introduction}

Transformational theory represents a challenging topic in contemporary ``mathemusical" research. It not only constitutes a turning point in the field of music analysis but also naturally leads to fundamental questions about the object-based vs. operation-based duality in the formalization of musical structures and processes. Starting from the original group-based contributions by David Lewin in the field of transformational music theory and analysis, the role of category theory has been quickly recognized \cite{Rahn2007,Mazzola2007,MazzolaTopos} as a useful foundation in order to build extended generalizations. Following previous work on the categorical formalization of Lewin's transformational theory and transformational networks in music analysis, this paper has two main goals. In a first part, we show how Lewin's definition of Generalized Interval Systems, even in a restricted form, implicitly defines diagrams of categories. In case all the condition of Lewin's definition are fullfilled, these diagrams involve category of elements, of which we discuss the importance both for further generalization and with respect to graph-based geometrical approach in music analysis. In a second part, we draw on these results to revisit the work on Klumpenhouwer networks and to provide a new definition of categorical transformational networks (``CT-Nets"), as well as that of morphisms between them.

\section{Revisiting Lewin's Generalized Interval Systems from a categorical perspective}

In this section, we revisit the original definition of Lewin's Generalized Interval Systems (GIS) \cite{Lewin1987} from a categorical perspective. We will show that Lewin's GIS inherently defines categories, which coincide with categories of elements for the group actions the GIS define due to the additional conditions imposed by Lewin.

\subsection{Lewin's Generalized Interval Systems}

We recall here Lewin's definition of Generalized Interval Systems (GIS) for transformational music analysis.

\begin{definition}[Lewin, 1987]
A Generalized Interval System (GIS) is a triple $(X,G,\text{int})$ where
\begin{itemize}
\item{$X$ is a set of musical elements,}
\vspace{0.2cm}
\item{$G$ is a group (the group of intervals for the GIS), and}
\vspace{0.2cm}
\item{$\text{int}$ is a function $X \times X \to G$}
\end{itemize}
\vspace{0.2cm}
such that
\vspace{0.2cm}
\begin{enumerate}
\item{\label{GIS-C1} for all $x$, $y$, and $z$ in $X$, $\text{int}(x,y)*\text{int}(y,z)=\text{int}(x,z)$, and}
\vspace{0.2cm}
\item{\label{GIS-C2}for all $x \in X$ and $g \in G$, there exists a unique $y \in X$ such that $\text{int}(x,y)=g$.}
\end{enumerate}
\end{definition}

We know from Vuza \cite{Vuza1988} and Kolman \cite{Kolman2004} that the data of a GIS is equivalent to the data of a simply transitive (right) group action of $G$ on the set $X$, or equivalently, of a representable functor $\mathbf{G} \to \mathbf{Sets}$ (where by an abuse of notation $\mathbf{G}$ is here the group $G$ considered as a single-object category). Note that if we replace condition (\ref{GIS-C1}) by $\text{int}(y,z)*\text{int}(x,y)=\text{int}(x,z)$, we obtain the data of a simply transitive left group action of $G$ on $X$.

The condition (\ref{GIS-C2}) is necessary for the existence of a group action. Indeed, if we keep only the condition (\ref{GIS-C1}), we may have an interval function $\text{int} \colon X \times X \to G$ which is not surjective. There would then be elements of $G$ for which we cannot define their action on $X$. However we will show in the following proposition that condition (\ref{GIS-C1}) alone in Lewin's definition actually defines the data of a category $\mathbf{C}$ along with a functor $\mathbf{C} \to \mathbf{G}$. Note that we use the alternative form of condition (\ref{GIS-C1}) to agree with the usual composition of morphisms in a category.

\begin{proposition}
\label{prop:lewin_category}
We consider a triple $(X,G,\text{int})$ where
\begin{itemize}
\item{$X$ is a set of musical elements,}
\vspace{0.2cm}
\item{$G$ is a group (the group of intervals for the GIS), and}
\vspace{0.2cm}
\item{$\text{int}$ is a function $X \times X \to G$ such that for all $x$, $y$, and $z$ in $X$, $\text{int}(y,z)*\text{int}(x,y)=\text{int}(x,z)$}
\end{itemize}
This defines the data of a category $\mathbf{C}$ along with a functor $\mathbf{C} \to \mathbf{G}$.
\end{proposition}
\begin{proof}
We construct $\mathbf{C}$ as the category having $X$ as its set of objects, and the set $\{\, (x,y) \colon x \to y \mid x \in X, y \in X\,\}$ as its set of morphisms. Composition of morphisms is straightforward, with $(y,z)(x,y)=(x,z)$. By the properties of the function $\text{int}$, we can extend it to a functor from $\mathbf{C}$ to $\mathbf{G}$.
\end{proof}

If condition \ref{GIS-C2} is added to the definition, then we obtain the following corollary.

\begin{definition}
Let $S \colon \mathbf{C} \to \mathbf{Sets}$ be a functor from a category $\mathbf{C}$ to $\mathbf{Sets}$. The category of elements $\int_\mathbf{C} S$ is defined as the category having
\begin{itemize}
\item{objects of the form $(X,x)$ with $X$ an object of $\mathbf{C}$ and $x$ an element of the set $S(X)$, and}
\item{morphisms between objects $(X,x)$ and $(Y,y)$ of the form $(x,f,y)$ with $f$ being a morphism of $\mathbf{C}$ such that $y=S(f)(x)$}
\end{itemize}
There is a canonical projection functor $\pi_S \colon \int_\mathbf{C} S \to \mathbf{C}$ sending each object $(X,x)$ to $X$.
\end{definition}

\begin{corollary}
Let $(X,G,\text{int})$ be a GIS. Then, the category $\mathbf{C}$ defined above is the category of elements $\int_\mathbf{G} S$ for the simply transitive group action $S \colon \mathbf{G} \to \mathbf{Sets}$
\end{corollary}
\begin{proof}
Immediate.
\end{proof}

In the following subsection, we explore how these results provide insights into the relationship between transformational music theory and category theory.

\subsection{From Lewin's Generalized Interval Systems to categories, and to graphs}

As we have seen above, the condition (\ref{GIS-C2}) enforces the existence of a simply transitive group action of $G$ on $X$, i.e. a representable functor $S \colon \mathbf{G} \to \mathbf{Sets}$. This gives rise to a diagram of categories of the form $\int_\mathbf{G} S \to \mathbf{G} \to \mathbf{Sets}$. Recent advances in transformational music theory have made use of various groups acting on sets of musical elements. Such groups do not always act in a simply transitive manner: this is the case for example of the $\text{T}/\text{I}$ group acting on pitch classes, or Hook's UTT group \cite{Hook2002} acting on major and minor triads. Although these do not constitute GIS, they nevertheless are useful and widely used tools for transformational analysis \cite{Lewin2006,Lewin2007}.

In this regard, dropping condition (\ref{GIS-C2}) altogether encourages us to shift our focus to categories rather than group actions, and to consider the category $\mathbf{C}$ and the associated functor $\mathbf{C} \to \mathbf{G}$ implied by Proposition \ref{prop:lewin_category} as a key notion. We will show in this section how we can generalize this observation by going beyond the definition of the interval function $\text{int}$.

As a first observation, it can be noted that for any group $G$ acting on a set $X$ we can always form the corresponding category of elements $\int_\mathbf{G} S$. While properly defining an interval function $\text{int}$ might be complicated if the group does not act in a simply transitive way, the definition of $\int_\mathbf{G} S$ and the associated functor $\int_\mathbf{G} S \to \mathbf{G}$ already contains the necessary properties of composability and associativity.

Recent research has emphasized he significance of extending beyond Lewin's original group-based approach, by, for instance, exploring the use of groupoids or categories in general \cite{Popoff2015,Popoff2019}. In this view, the group $\mathbf{G}$ and its action on a set is replaced with a general functor $S \colon \mathbf{C} \to \mathbf{Sets}$, and we can thus consider the corresponding category of elements $\int_\mathbf{C} S$ together with the canonical projection functor $\int_\mathbf{C} S \to \mathbf{C}$.

This generalization can be extended further by recognizing that the definition of the category of elements expands beyond functors to $\mathbf{Sets}$. In particular, this notion can be defined for functors $S \colon \mathbf{C} \to \mathbf{Rel}$, in which $\mathbf{Rel}$ is the category of finite sets and binary relations between them \cite{Sobocinski2012}. Such functors have been shown to occur in transformational music analysis, for example as in the algebraic formalization of the 'Cube Dance' graph \cite{Popoff2018,Popoff2022a}. One advantage of using the category of elements $\int_\mathbf{C} S$ is that it 'forgets' the target category of the functor $S$, whether it is $\mathbf{Rel}$ or $\mathbf{Sets}$, which allows one to consider them on an equal footing.

In passing, it can also be observed that the category of elements takes a role, in a discrete case, analog to that of the fundamental groupoid of a topological space: objects of the fundamental groupoid are points of the space, and 1-morphisms are paths between these points. This reconciles with Lewin's view of transformations as ways to pass from one object to another\footnote{During the redaction of this manuscript, the authors have been made aware by Dmitri Tymoczko of his current work on groupoids. Category of elements are encountered as a common point between his approach, rooted in geometrical and topological considerations, and ours, which stems from algebraic ones}.

\begin{figure}
\centering
\includegraphics[scale=1.3]{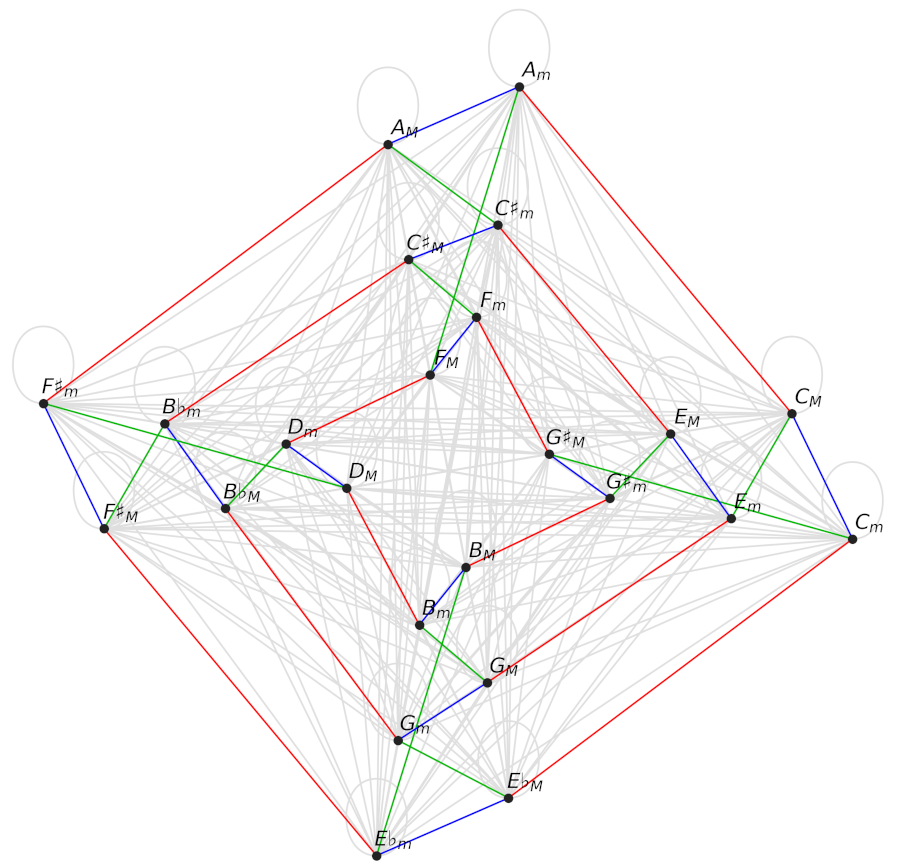}
\caption{The quiver corresponding to the category of elements for the simply transitive action of the $PLR$ group on the set of major and minor triads. Edges (corresponding to invertibles morphisms) are represented with undirected edges colored in gray for clarity, except for morphisms projecting on $R$ (red), $P$ (blue), and $L$ (green) in the $PLR$ group. This subquiver corresponds to the 'chicken-wire torus'. Best viewed in color.}
\label{fig:chicken_wire_graph_catelem}
\end{figure}

Finally, it can be noted that the notion of category of elements establishes a bridge between algebraic transformational approaches and geometrical graph-based approaches which have been studied extensively \cite{Douthett1998,Cohn2000,Cohn2012,Catanzaro2011,Piovan2013,Reenan2016,Rockwell2009}. Indeed, it is well known that there is a forgetful functor from $\mathbf{Cat}$ to $\mathbf{Quiv}$, the category of quivers, i.e. directed multigraphs ('multidigraphs'). The forgetful functor returns the underlying multidigraph of a category and forgets all information about composition of arrows. From this consideration, well-known structures in mathematical music theory can be retrieved. For example, the so-called 'chicken-wire torus' \cite{Douthett1998} is a sub-multidigraph of the multidigraph obtained by application of the forgetful functor $\mathbf{Cat} \to \mathbf{Quiv}$ on the category of elements $\int_{PLR} S$ for the usual action $S$ of the $PLR$ group on the set of major and minor triads. This is illustrated in Figure \ref{fig:chicken_wire_graph_catelem}, in which the multidigraph obtained from the category of elements is depicted in gray (since the action is simply transitive, the corresponding quiver is a complete one), with the 'chicken-wire torus' highlighted in color.

This colored multidigraph is a special case of the category-based general case, in which we consider the canonical projection functor $\int_\mathbf{C} S \to \mathbf{C}$ for a general functor $S \colon \mathbf{C} \to \mathbf{Sets}$. In this context, applying the forgetful functor $\mathbf{Cat} \to \mathbf{Quiv}$ induces additional structure, as shown in proposition below.

\begin{proposition}
Given a diagram $\mathbf{U} \to \mathbf{V}$, the forgetful functor $\mathbf{Cat} \to \mathbf{Quiv}$ returns a node- and edge-colored multidigraph corresponding to $\mathbf{U}$. The nodes are colored in the set of nodes of the multidigraph corresponding to $\mathbf{V}$, while the edges are colored in the set of edges of same multidigraph.
\end{proposition}
\begin{proof}
Immediate.
\end{proof}

In the above example, the category of elements for the usual action of the $PLR$ group on the set of major and minor triads gives an edge-colored multidigraph in which each edge is colored (labelled) by an element of the $PLR$ group. Note that we have colored the edges in Figure \ref{fig:chicken_wire_graph_catelem} in gray merely for clarity. The nodes have only one color since we are considering a single-object category.

\section{A new definition of transformational networks}

By following our shift of focus on categories of elements, we consider in this section a new definition of transformational networks which relates musical objects and transformations between them in a consistent diagram. Transformational networks first appeared in Lewin's work \cite{Lewin1987}, and so-called \textit{Klumpenhouwer networks} were developed by Klumpenhouwer \cite{Lewin1990,Klumpenhouwer1991,Klumpenhouwer1998} to specifically study pitch classes and their transformations by the $T/I$ group. In recent years, transformational networks have been revisited with new definitions from a categorical perspectives \cite{Rahn2007,Mazzola-Andreatta,Popoff2015,Popoff2016}. 

Our starting point is the definition of \cite{Popoff2015} which we recall below.

\begin{definition}
Let $\mathbf{C}$ be a category, and $S$ a functor from $\mathbf{C}$ to the category $\mathbf{Sets}$ of (small) sets. Let $\Delta$ be a small category and $R$ a functor from $\Delta$ to $\mathbf{Sets}$ with non-empy values. A {\it PK-net of form $R$ and of support $S$} is a 4-tuple $(R,S,F,\phi)$, in which $F$ is a functor from $\Delta$ to $\mathbf{C}$, and $\phi$ is a natural transformation from $R$ to $SF$.
\end{definition}

A typical PK-net is thus represented by the following diagram.

\begin{center}
	\begin{tikzpicture}
		\node (A) at (0,0) {$\Delta$};
		\node (B) at (3,0) {$\mathbf{C}$};
		\node (C) at (1.5,-2) {$\mathbf{Sets}$};
		\draw[->,>=latex] (A) -- (B) node[above,midway]{$F$};
		\draw[->,>=latex] (B) -- (C) node[right,midway] (D) {$S$} ;
		\draw[->,>=latex] (A) -- (C) node[left,midway]  (E) {$R$};
		\draw[->,>=latex,dashed] (A) to[out=-5,in=60,looseness=2.0] node[above,midway](F){} (C) ;
		\draw[-Implies,double distance=1.5pt,shorten >=8pt,shorten <=8pt] (E) to node[below,midway] {$\phi$} (F);
	\end{tikzpicture}
\end{center}

As per the discussion above, our objective is to shift the focus from $\mathbf{Sets}$-based functors to the corresponding categories of elements, which we do by leveraging the following theorem.

\begin{theorem}
\label{thm:catelements_square}
Assume we have two functors $R \colon \Delta \to \mathbf{Sets}$ and $S \colon \mathbf{C} \to \mathbf{Sets}$, as well as a functor $F \colon \Delta \to \mathbf{C}$ and a natural transformation $\phi \colon R \to SF$.
Then there exists a functor $G$ between the category of elements $\int_\Delta R$ and $\int_\mathbf{C} S$ such that the following diagram commute.

\begin{center}
	\begin{tikzpicture}
		\node (A) at (0,2) {$\int_\Delta R$};
		\node (B) at (2,2) {$\int_\mathbf{C} S$};
		\node (C) at (0,0) {$\Delta$};
		\node (D) at (2,0) {$\mathbf{C}$};

		\draw[->,>=latex] (A) -- (B) node[above,midway]{$H_F$};
		\draw[->,>=latex] (C) -- (D) node[above,midway]{$F$};
		\draw[->,>=latex] (A) -- (C) ;
		\draw[->,>=latex] (B) -- (D);
	\end{tikzpicture}
\end{center}

\end{theorem}
\begin{proof}
We construct the functor $H_F$ as follows.
\begin{itemize}
\item{For any object $(X,x)$ of $\int_\Delta R$, the functor $H_F$ sends it to $(F(X),\phi_X(x))$.}
\item{For any morphism $(x,f,y)$ between two objects $(X,x)$ and $(Y,y)$ of $\int_\Delta R$, the functor $H_F$ sends this morphism to $(\phi_X(x),F(f), \phi_Y(y))$.}
\end{itemize}
The construction in the first item is consistent, since by the definition of the natural transformation $\phi$, $\phi_X(x)$ is indeed an element of $SF(X)$. Similarly, the construction in the second item is also consistent, since for any $(x,y)$ such that $y=R(f)(x)$, we have $\phi_X(x)=SF(f)(\phi_Y(y))$.
We have $H_F((y,f',z)(x,f,y)) = H_F((y,f',z))H_F((x,f,y))$ by horizontal composition of the natural transformation.
\end{proof}

For a given PK-net, we thus have the following diagram.

\begin{center}
	\begin{tikzpicture}
		\node (A) at (0,	0) {$\Delta$};
		\node (B) at (3,0) {$\mathbf{C}$};
		\node (C) at (1.5,-2) {$\mathbf{Sets}$};

		\node (A2) at (0,2) {$\int_\Delta R$};
		\node (B2) at (3,2) {$\int_\mathbf{C} S$};

		\draw[->,>=latex] (A2) -- (B2) node[above,midway]{$H_F$};
		\draw[->,>=latex] (A2) -- (A) ;
		\draw[->,>=latex] (B2) -- (B);

		\draw[->,>=latex] (A) -- (B) node[above,midway]{$F$};
		\draw[->,>=latex] (B) -- (C) node[right,midway] (D) {$S$} ;
		\draw[->,>=latex] (A) -- (C) node[left,midway]  (E) {$R$};
		\draw[->,>=latex,dashed] (A) to[out=-5,in=60,looseness=2.0] node[above,midway](F){} (C) ;
		\draw[-Implies,double distance=1.5pt,shorten >=8pt,shorten <=8pt] (E) to node[below,midway] {$\phi$} (F);
	\end{tikzpicture}
\end{center}

The definition of PK-nets allows us to define networks in which nodes are sets, and not necessarily singletons. In this diagram, the role of the functor $\int_\Delta R \to \Delta$ is to remember which objects of $\int_\Delta R$ belong to the same set. The fact that the top square is commutative ensures that objects corresponding to elements of the same set are transformed by the same transformations.

While this possibility gives more generality to transformational networks, the majority of actual transformational analyses in the musicological literature uses networks with only one musical object per node, i.e. singletons. In this case, $\int_\Delta R$ is the same category as $\Delta$, and since the data of the content of sets is no longer relevant, this makes the functor $\int_\Delta R \to \Delta$ somewhat unnecessary. We will thus now consider the specific case of singleton-based networks and drop this functor altogether. In the same movement, we can rename $\int_\Delta R$ as $\Delta$ to indicate its new role as a category defining the skeleton of the transformational network to be considered. We therefore obtain a diagram of the following form.

\begin{center}
	\begin{tikzpicture}
		\node (A) at (0,0) {$\Delta$};
		\node (B) at (2,0) {$\int_\mathbf{C} S$};
		\node (C) at (4,0) {$\mathbf{C}$};

		\draw[->,>=latex] (A) -- (B) node[above,midway]{$H$};
		\draw[->,>=latex] (B) -- (C) node[above,midway] (E) {$\pi_S$} ;
	\end{tikzpicture}
\end{center}

This diagram encodes transformational networks with a single musical element on each node. The nodes and arrows of the network are given by the objects and morphisms of the category $\Delta$. The functor $H$ maps the objects of $\Delta$ to the objects of $\int_\mathbf{C} S$ which, as we have seen above, corresponds to musical elements. In other words, the functor $H$ labels the nodes of the network with musical elements. In addition, it maps the arrows of the network to morphisms of $\int_\mathbf{C} S$, which are then given a label in $\mathbf{C}$ (for example, group elements if $\mathbf{C}$ is a group) by the functor $\pi_S \colon \int_\mathbf{C} S \to \mathbf{C}$. Thus, the functor $\pi_S \circ H$ labels the arrows of $\Delta$ with morphisms (transformations) of $\mathbf{C}$.

In a final movement, we can observe that, since we do not consider $\mathbf{Sets}$-based (or $\mathbf{Rel}$-based) functors anymore, we can consider all categories in our definition of a transformational network, and not just category of elements. We thus obtain the following definition.

\begin{definition}
A categorical transformational network (CT-Net) is a diagram of the form
\begin{center}
	\begin{tikzpicture}
		\node (A) at (0,0) {$\Delta$};
		\node (B) at (2,0) {$\mathbf{C}_{\text{el}}$};
		\node (C) at (4,0) {$\mathbf{C}_{T}$};

		\draw[->,>=latex] (A) -- (B) node[above,midway]{$H$};
		\draw[->,>=latex] (B) -- (C) node[above,midway] (E) {$\pi$} ;
	\end{tikzpicture}
\end{center}
where
\begin{itemize}
\item{$\Delta$ is a category representing the skeleton of the network,}
\item{$\mathbf{C}_{\text{el}}$ is a category whose objects represent musical objects of interest, with morphisms between them, and}
\item{$\mathbf{C}_{T}$ is a category whose objects represent classes of musical objects, and musical transformations between them.}
\end{itemize}
The nodes of the network are labelled by the images of the objects of $\Delta$ by $H$, and the edges of the network are labelled by the images of the morphisms of $\Delta$ by $\pi \circ H$.
\end{definition}

Note that we do not \textit{a priori} impose conditions on the categories $\mathbf{C}_{\text{el}}$ and $\mathbf{C}_{T}$, nor on the functor $\pi$. In fact, $\mathbf{C}_{T}$ may well be $\mathbf{C}_{\text{el}}$ itself, in which case transformations are individualized for each musical object.

\subsubsection{Morphisms of CT-Nets}

In this section, we revisit the notions of morphisms of transformational networks and explore their implications for diagrams involving categories of elements. We will build upon the morphisms of PK-Nets defined in previous work, and in this view we will revert to category of elements for functors $\mathbf{C} \to \mathbf{Sets}$ before giving the final definition in the general case.

We first consider the notion of \textit{complete homographies} introduced in the formalization of PK-Nets. We recall their definition below.

\begin{definition}
Let $(R,S,F,\phi)$ be a PK-Net, and assume we have a functor $S' \colon \mathbf{C'} \to \mathbf{Sets}$. A complete homography is a pair $(N,\nu)$ with
\begin{itemize}
\item{$N \colon \mathbf{C} \to \mathbf{C'}$, and}
\item{$\nu \colon SF \to S'FN$ such that $\nu = \tilde{\nu}F$, where $\tilde{\nu}$ is a natural transformation from $S$ to $S'N$}
\end{itemize}
Applying $(N,\nu)$ on the PK-net $(R,S,F,\phi)$ transforms it into the new PK-Net $(R,S',FN,\nu \circ \phi)$.
\end{definition}

When $\mathbf{C}$ and $\mathbf{C'}$ are groups, with simply transitive group actions, this is known as a morphism of GIS \cite{Fiore2013}. For a given complete homography, we can apply Theorem \ref{thm:catelements_square} on the diagram formed by $S$ and $S'N$. This defines a functor $J_N$ from $\int_\mathbf{C} S$ to $\int_\mathbf{C'} S'$. Starting from a transformational network
\begin{center}
	\begin{tikzpicture}
		\node (A) at (0,0) {$\Delta$};
		\node (B) at (2,0) {$\int_\mathbf{C} S$};
		\node (C) at (4,0) {$\mathbf{C}$};

		\draw[->,>=latex] (A) -- (B) node[above,midway]{$H$};
		\draw[->,>=latex] (B) -- (C) node[above,midway] (E) {$\pi_S$} ;
	\end{tikzpicture}
\end{center}
we thus obtain a new transformational network, as shown in the diagram below.

\begin{center}
	\begin{tikzpicture}
		\node (A) at (0,0) {$\Delta$};
		\node (B) at (2,0) {$\int_\mathbf{C} S$};
		\node (C) at (4,0) {$\mathbf{C}$};

		\node (B2) at (2,-2) {$\int_\mathbf{C'} S'$};
		\node (C2) at (4,-2) {$\mathbf{C'}$};

		\draw[->,>=latex] (A) -- (B) node[above,midway]{$H$};
		\draw[->,>=latex] (B) -- (C) node[above,midway] (E) {$\pi_S$} ;

		\draw[->,>=latex] (A) -- (B2) node[left,midway]  (G) {$J_N \circ H$};
		\draw[->,>=latex] (B) -- (B2) node[left,midway]  (H) {$J_N$};
		\draw[->,>=latex] (B2) -- (C2) node[below,midway] (E) {$\pi_{S'}$} ;
		\draw[->,>=latex] (C) -- (C2) node[left,midway]  (I) {$N$};
		

	\end{tikzpicture}
\end{center}

Secondly, we consider \textit{local homographies}, which are PK-net homographies $(N,\nu)$ between PK-nets sharing the same functors $R$ and $S$ such that the natural transformation $\nu \colon SF \to SF'$ can be expressed as $\nu = S\hat{\nu}$, where $\hat{\nu}$ is a natural transformation from $F$ to $F'=NF$. As seen in a recent paper \cite{Popoff2019}, we can even drop the condition that $F'=NF$ and consider natural transformations $\hat{\nu} \colon F \to F'$, which are extended by $S$ such that $\phi' = (S\hat{\nu})\phi$. In such a case, it can readily be seen that the natural transformation $\hat{\nu} \colon F \to F'$ induces a natural transformation $\eta_{\hat{\nu}} \colon H_F \to H_{F'}$ between the corresponding functors $H_F \colon \int_\Delta R \to  \int_\mathbf{C} S$ and $H_{F'} \colon \int_\Delta R \to  \int_\mathbf{C} S$. We then obtain the following diagram.

\begin{center}
	\begin{tikzpicture}
		\node (A) at (0,0) {$\Delta$};
		\node (B) at (2,0) {$\int_\mathbf{C} S$};
		\node (C) at (4,0) {$\mathbf{C}$};

		\draw[->,>=latex] (A) to[bend left=20] node[above,midway] (U) {$H_F$} (B) ;
		\draw[->,>=latex] (A) to[bend right=20] node[below,midway] (U2) {$H_{F'}$} (B) ;
		\draw[->,>=latex] (B) -- (C) node[above,midway] (E) {$\pi_S$} ;
		\draw[-Implies,double distance=1.5pt,shorten >=2pt,shorten <=2pt] (U) to node[right,midway] {$\eta_{\hat{\nu}}$} (U2);
	\end{tikzpicture}
\end{center}

By combining these two notions, and by generalizing it through the possibility of changing the category $\Delta$, we propose a unified notion of morphism of CT-Nets as defined below.

\begin{definition}
A morphism between two categorical transformational networks
\begin{center}
	\begin{tikzpicture}
		\node (A) at (0,0) {$\Delta$};
		\node (B) at (2,0) {$\mathbf{C}_{\text{el}}$};
		\node (C) at (4,0) {$\mathbf{C}_{T}$};
		
		\draw[->,>=latex] (A) -- (B) node[above,midway]{$H$};
		\draw[->,>=latex] (B) -- (C) node[above,midway] (E) {$\pi$} ;
	\end{tikzpicture}
\end{center}
and
\begin{center}
	\begin{tikzpicture}
		\node (A) at (0,0) {$\Delta'$};
		\node (B) at (2,0) {$\mathbf{C'}_{\text{el}}$};
		\node (C) at (4,0) {$\mathbf{C'}_{T}$};

		\draw[->,>=latex] (A) -- (B) node[above,midway]{$H'$};
		\draw[->,>=latex] (B) -- (C) node[above,midway] (E) {$\pi'$} ;
	\end{tikzpicture}
\end{center}
is defined as a 4-tuple $(I,N,J_N,\nu)$ such that we have the following commutative diagram.
\begin{center}
	\begin{tikzpicture}
		\node (A) at (0,0) {$\Delta$};
		\node (B) at (2.5,0) {$\mathbf{C}_{\text{el}}$};
		\node (C) at (5,0) {$\mathbf{C}_{T}$};

		\draw[->,>=latex] (A) -- (B) node[above,midway]{$H$};
		\draw[->,>=latex] (B) -- (C) node[above,midway] (E) {$\pi$} ;

		\node (A2) at (0,-2.5) {$\Delta'$};
		\node (B2) at (2.5,-2.5) {$\mathbf{C'}_{\text{el}}$};
		\node (C2) at (5,-2.5) {$\mathbf{C'}_{T}$};

		\draw[->,>=latex] (A2) -- (B2) node[below,midway]{$H'$};
		\draw[->,>=latex] (B2) -- (C2) node[below,midway] (E) {$\pi'$} ;
		
		\draw[->,>=latex] (A) -- (A2) node[left,midway]{$I$};
		\draw[->,>=latex] (B) -- (B2) node[right,midway]{$J_N$};
		\draw[->,>=latex] (C) -- (C2) node[right,midway]{$N$};
		
		\draw[->,>=latex,dashed] (A) to[out=-5,in=95,looseness=0.9] node[above,midway](U){} (B2) ;
		\draw[->,>=latex,dashed] (A) to[out=-85,in=175,looseness=0.9] node[above,midway](V){} (B2) ;
		\draw[-Implies,double distance=1.5pt,shorten >=2pt,shorten <=7pt] (U) to node[below,midway] {$\nu$} (V);

	\end{tikzpicture}
\end{center}	
Morphisms of CT-Nets are composable, via the appropriate stacking of the corresponding diagrams.
\end{definition}

In the above definition of a morphism of CT-Nets, the rightmost square represents a global transformation of the musical elements and the corresponding transformations between them, independently of the specific network $\Delta$. As detailed above, if $\mathbf{C}_{\text{el}}$ is the category of elements for a given group action, this corresponds to a morphism of a group action. The leftmost square represents a local transformation of the network nodes: the components of the natural transformation $\nu$ are morphisms in $\mathbf{C'}_{\text{el}}$, which can be given labels via $\pi'$ in $\mathbf{C'}_{T}$. We give below an example of such a combination of global and local transformation.

\begin{example}
The two transformational networks represented below are $\langle T_2\rangle$-isographic in the sense of Klumpenhouwer.
\begin{center}
\begin{tikzpicture}
	\node (E) at (4,0) {$B$};
	\node (F) at (6,0) {$C_\sharp$};
	\node (G) at (4,-2) {$E$};
	\node (H) at (6,-2) {$G_\sharp$};
	\draw[->,>=latex] (E) -- (F) node[above,midway]{$T_2$};
	\draw[->,>=latex] (F) -- (H) node[right,midway]{$I_9$};
	\draw[->,>=latex] (E) -- (G) node[left,midway]{$I_3$};
	\draw[->,>=latex] (G) -- (H) node[below,midway]{$T_4$};

	\node (I) at (0,0) {$G_\sharp$};
	\node (J) at (2,0) {$B_\flat$};
	\node (K) at (0,-2) {$F$};
	\node (L) at (2,-2) {$A$};
	\draw[->,>=latex] (I) -- (J) node[above,midway]{$T_2$};
	\draw[->,>=latex] (J) -- (L) node[right,midway]{$I_9$};
	\draw[->,>=latex] (I) -- (K) node[left,midway]{$I_3$};
	\draw[->,>=latex] (K) -- (L) node[below,midway]{$T_4$};
	\end{tikzpicture}
\end{center}

However, no morphism of the $T/I$ group action can transform the first into the second. This can however be achieved with a combination of a global transformation with a local one, as shown below.

\begin{center}
\begin{tikzpicture}
	\node (A) at (-5,0) {$G_\sharp$};
	\node (B) at (-3,0) {$B_\flat$};
	\node (C) at (-5,-2) {$F$};
	\node (D) at (-3,-2) {$A$};
	\draw[->,>=latex] (A) -- (B) node[above,midway]{$T_2$};
	\draw[->,>=latex] (B) -- (D) node[right,midway]{$I_7$};
	\draw[->,>=latex] (A) -- (C) node[left,midway]{$I_1$};
	\draw[->,>=latex] (C) -- (D) node[below,midway]{$T_4$};

	\node (E) at (5,0) {$B$};
	\node (F) at (7,0) {$C_\sharp$};
	\node (G) at (5,-2) {$E$};
	\node (H) at (7,-2) {$G_\sharp$};
	\draw[->,>=latex] (E) -- (F) node[above,midway]{$T_2$};
	\draw[->,>=latex] (F) -- (H) node[right,midway]{$I_9$};
	\draw[->,>=latex] (E) -- (G) node[left,midway]{$I_3$};
	\draw[->,>=latex] (G) -- (H) node[below,midway]{$T_4$};

	\node (I) at (0,0) {$E_\flat$};
	\node (J) at (2,0) {$F$};
	\node (K) at (0,-2) {$C$};
	\node (L) at (2,-2) {$E$};
	\draw[->,>=latex] (I) -- (J) node[above,midway]{$T_2$};
	\draw[->,>=latex] (J) -- (L) node[right,midway]{$I_9$};
	\draw[->,>=latex] (I) -- (K) node[left,midway]{$I_3$};
	\draw[->,>=latex] (K) -- (L) node[below,midway]{$T_4$};

	\node (M) at (2.75,-1) {};
	\node (N) at (4.25,-1) {};
	\node (O) at (-2.25,-1) {};
	\node (P) at (-0.75,-1) {};
	\draw[->,>=latex] (M) -- (N);
	\draw[->,>=latex] (O) -- (P) node[above,midway]{$\langle T_2 \rangle$};

	\draw[->,>=latex,dashed,red] (I) to[bend left] node[above,midway]{\color{red}{$T_8$}} (E) ;
	\draw[->,>=latex,dashed,red] (J) to[bend left] node[above,midway]{\color{red}{$T_8$}} (F) ;
	\draw[->,>=latex,dashed,red] (K) to[bend right] node[below,midway]{\color{red}{$T_4$}} (G) ;
	\draw[->,>=latex,dashed,red] (L) to[bend right] node[below,midway]{\color{red}{$T_4$}} (H) ;
\end{tikzpicture}
\end{center}
\end{example}

\begin{example}
As is well-known in category theory, the functor category $\mathbf{Sets}^{\mathbf{C}}$ is a topos. Therefore, if one works with functors $S \colon \mathbf{C} \to \mathbf{Sets}$, with $\mathbf{C}_{\text{el}} = \int_\mathbf{C} S$ and $\mathbf{C}_{T} = \mathbf{C}$, characteristic morphisms readily provide global transformations \cite{Popoff2015}.

Assume $A \subset S$ is a subobject of $S$. We then have a characteristic morphism $\chi_A \colon S \to \Omega$ from $S$ to $\Omega$, the subobject classifier of $\mathbf{Sets}^{\mathbf{C}}$. This translates into a global transformation of the corresponding category of elements, and to a morphism of CT-Nets as shown in the diagram below.

\begin{center}
	\begin{tikzpicture}
		\node (A) at (0,0) {$\Delta$};
		\node (B) at (2,0) {$\int_\mathbf{C} S$};
		\node (C) at (4,0) {$\mathbf{C}$};

		\node (B2) at (2,2) {$\int_\mathbf{C} A$};

		\node (B3) at (2,-2) {$\int_\mathbf{C} \Omega$};
		\node (C3) at (4,-2) {$\mathbf{C}$};

		\draw[->,>=latex] (A) -- (B) node[above,midway]{$H$};
		\draw[->,>=latex] (B) -- (C) node[above,midway] (E) {$\pi_S$} ;

		\draw[->,>=latex] (B2) -- (B) node[left,midway]  (H) {$J$};
		\draw[->,>=latex] (C) -- (C3) node[right,midway]  (I) {$\text{id}$};
		\draw[->,>=latex] (B) -- (B3) node[left,midway]  (J) {$J_{\chi_A}$};
		\draw[->,>=latex] (B3) -- (C3) node[below,midway]  (J) {$\pi_\Omega$};
		
	\end{tikzpicture}
\end{center}
A new transformational network
\begin{center}
	\begin{tikzpicture}
		\node (A) at (0,0) {$\Delta$};
		\node (B) at (2,0) {$\int_\mathbf{C} \Omega$};
		\node (C) at (4,0) {$\mathbf{C}$};

		\draw[->,>=latex] (A) -- (B) node[above,midway]{$J_{\chi_A} \circ H$};
		\draw[->,>=latex] (B) -- (C) node[above,midway] (E) {$\pi_\Omega$} ;
	\end{tikzpicture}
\end{center}	
is thus obtained.
\end{example}

Since our definition of transformational networks is entirely categorical, all constructions in $\mathbf{Cat}$ can in fact be used to construct new CT-Nets. As shown above, pullbacks and pushforwards readily provide morphisms of CT-Nets. Furthermore, as mentionned above, the passage from $\mathbf{Cat}$ to $\mathbf{Quiv}$ is functorial which allows to translate these constructions into the corresponding multidigraphs (and morphisms between them).

\section{Conclusions}

In this work, we have revisited Lewin's original framework of Generalized Interval Systems from a categorical perspective, by stressing out the importance of diagrams $\mathbf{C}_{\text{el}} \to \mathbf{C}_{T}$, in which $\mathbf{C}_{\text{el}}$ is a category of musical objects which maps onto $\mathbf{C}_{T}$, a category of musical transformations. We have shown that Lewin's definition of GIS implictly contains the definition of such a diagram, and that the condition for simple transitivity turns $\mathbf{C}_{\text{el}}$ into the category of elements $\int_\mathbf{G} S$ for the corresponding group action $S \colon \mathbf{G} \to \mathbf{Sets}$. This consideration extends beyond group and group actions however, and category of elements $\int_\mathbf{C} S$ for a given functor $S \colon \mathbf{C} \to \mathbf{Sets}$ or even $S \colon \mathbf{C} \to \mathbf{Rel}$ may be considered. Ultimately, one may consider diagrams $\mathbf{C}_{\text{el}} \to \mathbf{C}_{T}$ in which $\mathbf{C}_{\text{el}}$ is not necessarily a category of elements.

Drawing from previous work on the formalization of transformational networks, we have proposed a new categorical definition of transformational networks, hence called CT-Nets, as diagrams $\Delta \to \mathbf{C}_{\text{el}} \to \mathbf{C}_{T}$, in which $\Delta$ represents the skeleton of the network. The definition operates directly in $\mathbf{Cat}$, and thus unifies the consideration of $\mathbf{Sets}$- and $\mathbf{Rel}$-based functors at once. The definition of morphisms of CT-Nets results from previous considerations of global and local morphisms, the former corresponding to global transformations of the musical action, independently of the network being considered, the latter corresponding to individual node transformations. As diagrams, CT-Nets can undergo the machinery of all constructions in $\mathbf{Cat}$ (pullbacks, pushforwards, equalizers, etc.) to give new networks. In addition, the forgetful functor $\mathbf{Cat} \to \mathbf{Quiv}$ readily gives the corresponding multidigraphs (and morphisms between them), bridging algebraic and geometrical graph-based methods in music theory.

The new perspective on transformational networks introduced in this paper also provides lead for future work. For example, while category of elements can be considered for functors $\mathbf{C} \to \mathbf{Rel}$, Theorem \ref{thm:catelements_square} cannot be applied to relational PK-nets as introduced in \cite{Popoff2018} if the natural transformation $\phi$ is a binary relation and not simply a function. In such a case, the machinery of profunctors \cite{Benabou1973} may be used between the corresponding category of elements, with new implications regarding the musicological meaning of transformations between musical objects.

\section*{Acknowledgement}
We thank Andr\'ee Ehresmann for many fruitful discussions on this topic.

\bibliography{arxiv_hidden_categories}
\bibliographystyle{alpha}
\addcontentsline{toc}{section}{References}

\end{document}